\documentclass[11pt,reqno,a4paper]{amsart}
\usepackage{amsmath,amsfonts,amssymb,amsthm,amsxtra}
\usepackage{tikz-cd}
\usepackage{tikz}
\usetikzlibrary{positioning,chains,fit,shapes,calc,arrows,patterns,external,shapes.callouts,graphs}

\newtheorem{theorem}{Theorem}[section]
\newtheorem{lemma}{Lemma}[section]
\newtheorem{remark}{Remark}[section]
\newtheorem{proposition}{Proposition}[section]
\newtheorem{example}{Example}[section]
\newtheorem{definition}{Definition}[section]
\newtheorem{corollary}{Corollary}[section]

\newcommand{\bo}{[\![}
\newcommand{\bc}{]\!]}

\newcommand{\extd}{d_L}
\newcommand{\extdd}{d_{L^*}}
\newcommand{\Lie}{\mathcal{L}}

\newcommand{\llangle}{\langle \!\! \langle}
\newcommand{\rrangle}{\rangle \!\! \rangle}

\begin{document}
\title{T-Duality and Lie bi-algebroid structures}
\author{Alexander Cardona}
\author{Juan Jos\'e Villamar\'in}
\address{Mathematics Department
\\Universidad de Los Andes
\\Carrera 1 No. 18A-10
\\ 111711 Bogot\'a, Colombia.}
\email{acardona@uniandes.edu.co}
\email{jj.villamarin@uniandes.edu.co}
\date{July 26, 2022}

\maketitle

\begin{abstract}
We show that the geometric notion of duality behind $T$-duality, between two string theories on different manifolds $E, \hat{E}$ in the sense of \cite{BHM1}\cite{BHM2}, is precisely that of Lie bialgebroids due to Mackenzie and Xu \cite{MX}. 
\end{abstract}

\vspace{0.5cm}

\textbf{MSC(2000)}: 53C15, 53D17.

\textbf{Keywords}: Lie algebroids, Lie bialgebroids,  T-duality, Poisson brackets.

\vspace{0.5cm}

\section{Introduction}
Lie bialgebroids are infinitesimal counterparts of Poisson groupoids, they were introduced by Kirill Mackenzie and Ping Xu in \cite{MX} and generalize both the double structure of Poisson manifolds (i.e. pairs $( T M , T ^ { * } M )$, where $M$ is a Poisson manifold) and that of Lie bialgebras (in the sense of V.G. Drinfel'd, see \cite{V}). Lie bialgebroids encode the most natural definition of \emph{duality} for Lie algebroids, these are pairs of dual vector bundles over a base manifold $M$ defining Lie algebroids $(A, A^{*})$ satisfying, in addition, the compatibility condition that the coboundary operator $d _{A}*: A \to \bigwedge\nolimits^2 A$, associated to $A^*$, be a derivation of the graded Lie bracket $[ \cdot , \cdot ]_A$ on $\Gamma(A)$. The main goal of this paper is to show that the geometric notion of duality behind the so-called \emph{topological $T$-duality} between two string theories on different manifolds $E, \hat{E}$, in the sense of \cite{BHM1} \cite{BHM2}, is precisely that of Lie bialgebroids. In this setup the manifolds $E$ and $\hat{E}$ are principal $T$-bundles, where $T = \mathbb{S}^1 \times \cdots \times \mathbb{S}^1$ is a torus, over the same base manifold $M$.  \\
\\
$T$-duality is an important tool in string theory. Among other features, it gives rise to an equivalence between string theories which in their low-energy limit might look different, but exhibit a one-to-one correspondence between fields, states and other defining ingredients of the corresponding theories. In the context of principal bundles, the smooth version of $T$-duality is a duality of torus bundles induced by closed invariant 3-forms as follows \cite{BHM1} \cite{BHM2} (see also \cite{CG}).  Let us denote by $E$ and $\hat{E}$ two principal $T$-bundles, where $T$ denotes a $k$-torus with Lie algebra $\mathfrak{t}$, over the same base manifold $M$. The source of the duality will be a non-degenerate pairing $\Phi_o: \mathfrak{t} \times \mathfrak{t}^* \rightarrow \mathbb{R}$, where $\mathfrak{t}^*$ denotes the dual of the Lie algebra $\mathfrak{t}$, and it will be realized on the so-called \emph{correspondence space} $E\times_M \hat{E}$ (the fiber product) associated to the space-time manifolds $E$ and $\hat{E}$. Together with the principal $T$-bundles $E \overset{\pi}{\rightarrow} M$ and $\hat{E} \overset{\hat{\pi}}{\rightarrow} M$ we will consider closed $T$-invariant 3-forms $H$ and $\hat{H}$, with integral periods on $E$ and $\hat{E}$, respectively. The correspondence space define principal fibrations over $E$ and $\hat{E}$, and the pull-back of the $T$-invariant 3-forms $H$ and $\hat{H}$ gives rise to an element $\Phi \in \Omega^2(E\times_M \hat{E})$ which plays a double role in the construction: geometrically it will define the double structure encoding the $T$-duality, topologically it intertwines between the twisted cohomologies (and twisted K-theories) of $(E, H),$ and $(\hat{E}, \hat{H})$ (see, e.g. \cite{BHM1} \cite{BHM2}\cite{BS}).
\begin{equation}
\begin{tikzcd}[row sep=2.5em,column sep={4em,between origins},nodes in empty cells]
&
(E\times_M \hat{E}, p^* H- \hat{p}^*\hat{H} = d \Phi)
  \arrow[ld,swap,"p"]
  \arrow[rd,"\hat{p}"] &
\\
(E, H)
  \arrow[rd,swap,"\pi"] &  &
(\hat{E}, \hat{H}) 
  \arrow[ld,"\hat{\pi}"] 
\\
& M. &
\end{tikzcd}
\end{equation}

\begin{definition}
We say that the pairs $(E, H),$ and $(\hat{E}, \hat{H})$ are $T$-dual if 
    \begin{enumerate}
        \item $p^* H- \hat{p}^*\hat{H}  = d\Phi$ in the correspondence space $E\times_M \hat{E}$, for a $T \times \hat{T} \cong T^2$-invariant form $\Phi$.
        \item The form $\Phi$ restricted to the fibers of $E\times_M\hat{E}$ gives rise to the non-degenerate pairing $\Phi_o: \mathfrak{t} \times \mathfrak{t}^* \rightarrow \mathbb{R}$.
     \end{enumerate}
\end{definition}

Given a principal $T$-bundle $E \overset{\pi}{\rightarrow} M$ and a closed invariant 3-form $H$ on $E$, with integral periods, a pair $(E, H)$ is called \emph{$T$-dualizable} if there exists a closed 2-form $\hat{F} \in \Omega^2(M, \mathfrak{t}^*)$ with integral periods such that,  for all $X \in \mathfrak{t}$,
    \begin{equation}\label{eq:FCurv}
       i_X H = \hat{F}(X),
    \end{equation}
where $\hat{F}(X)$ means the natural pairing (with values in $\Omega^2(M)$) between $\mathfrak{t}$ and $\Omega^2(M,\mathfrak{t}^*)$. It has been shown that, given a $T$-dualizable pair
$(E, H)$, then there exists a $T$-dualizable pair $(\hat{E}, \hat{H})$ that is $T$-dual to $(E,H)$. As a matter of fact, in the case of circle bundles, the 2-forms defined by (\ref{eq:FCurv}) corresponding to $T$-dual pairs have a natural geometric interpretation as curvatures of the principal $\mathbb{S}^1$-bundles $E$ and $\hat{E}$. Namely,  there exist curvature forms $F$ and $\hat{F}$ on $E$ and $\hat{E}$, respectively, such that at the level of cohomology
    \begin{equation}\label{eq:2}
        \int_{\mathbb{S}^1} H = \hat{F} \hspace{8pt} \text{and} \hspace{8pt} \int_{\hat{\mathbb{S}}^1} \hat{H} = F.  
    \end{equation}
In that case $\Phi = - A \wedge \hat{A}$, where $A$ and $\hat{A}$ denote connections on the $\mathbb{S}^1$-bundles $E$ and $\hat{E}$ satisfying $F=dA$ and $\hat{F}=d \hat{A}$, respectively (see \cite{BHM1} \cite{BHM2} and references therein).\\
\\
Given $T$-dual pairs $(E, H)$ and $(\hat{E}, \hat{H})$, we will denote by $A({E})$ and $A(\hat{E})$, respectively, their associated Atiyah algebroids. Recall that, given a principal $G$-bundle $P \overset{\pi}{\rightarrow} M$, the sub-bundle induced by the right-invariant vector fields on $P$, denoted by $\mathfrak{X}_R(P)$, is integrable and can be identified with the space of sections of the bundle $TP/G$. Since the projection map $TP \rightarrow P$ is $G$-equivariant we obtain a map $TP/G \rightarrow P/G \simeq M$ so that we can see $TP/G$ as a vector bundle over $M$. If we consider $\mathrm{a} = d\pi : TP/G \rightarrow TM$ and the bracket as the restriction of the Lie bracket from $TP$, we induce a Lie algebroid structure on $TP/G$, called the Atiyah algebroid of $P$, and denoted by $A(P)$. Thus, we obtain the following short exact sequence of vector bundles over $M$
$$
                0 \rightarrow Ad(P) \rightarrow A(P) \rightarrow TM \rightarrow 0,
$$
splittings of which correspond with connections on $P$ (see \cite{A}\cite{M} and references therein). Given any Lie algebroid $A$, a differential $d_{ A }$ is defined on the graded algebra of sections of the exterior algebra of the dual vector bundle, $\Gamma ( \bigwedge A^{ * } )$, called the Lie algebroid differential of $A$. This differential generalizes the Chevalley–Eilenberg cohomology operator on $\wedge  \mathfrak {g}^{*} $ (when $A$ is a Lie algebra $\frak g$), the usual de Rham differential on forms (when $A = T M$) and the Lichnerowicz–Poisson differential $[\Pi , \cdot ]$ on multi-vector fields on $M$ (when $A = T^{*} M$ is the cotangent bundle of a Poisson manifold $(M, \Pi)$), among others.
\\
\\
The total spaces of circle bundles in $T$-dual pairs $(E, H)$ and $(\hat{E}, \hat{H})$ are contact manifolds, and both $F$ and $\hat{F}$ define symplectic forms on the base space $M$ (see \cite{BW}), we will use the associated Poisson structure to prove that the Atiyah algebroids $L= A({E})$ and $\hat{L}= A(\hat{E})$ over $M$ define a Lie bi-algebroid, i.e. they are dual 
(with anchor $a$ and bracket $[\cdot, \cdot]$, $a_*$ and bracket $[\cdot, \cdot]_*$, respectively) and $d_L$ is a derivation of the Schouten bracket of $\hat{L}$, so that for all $\alpha, \beta \in \Gamma(\hat{L})$:
\begin{equation}
            \extd[\alpha, \beta]_* = [\extd \alpha, \beta]_* + [\alpha, \extd\beta]_*.
\end{equation}
The duality pairing giving the identification of  $\hat{L}$ with $L^*$ is the defined by 
\begin{equation}\label{Eq:Duality}
 \llangle X, Y \rrangle  :=\Phi(\Tilde{X}^v,\Tilde{Y}^v) + F(X,Y),
\end{equation}
where $\Tilde{X}$ and $\Tilde{Y}$ denote horizontal lifts to the correspondence space $E\times_M\hat{E}$ of the vector fields on each fibration. \\ \\
In section \ref{S:Lie bialgebroids and Poisson algebras} we recall the basic facts and definitions on Lie bialgebroids, and the way they can be used to define Poisson algebras of functions on smooth manifolds. In section \ref{T-duality} we show that $T$-duality of circle bundles can be naturally interpreted as a Lie bialgebroid induced by the Atiyah algebroids associated to the $T$-dual pairs, and we will address some of the implications of this results.

\section{Lie bialgebroids and Poisson algebras}\label{S:Lie bialgebroids and Poisson algebras}

Consider a Lie algebroid $A$ over a smooth manifold $M$, with anchor $\mathrm{a}: \Gamma(A) \to \mathfrak{X}(M)$ and bracket $[\cdot, \cdot]$, whose dual $A^*$ is also equipped with a Lie algebroid structure, say $\mathrm{a}^*: \Gamma(A^*) \to \mathfrak{X}(M)$ and bracket $[\cdot, \cdot]_*$: 
$$
\begin{tikzcd}
A  \arrow[r, "\mathrm{a}"] \arrow[rd] &  TM  \arrow[d]  & A^* \arrow[l,swap, "\mathrm{a}^*"] \arrow[dl]&  \\
& M.  
\end{tikzcd}
$$
If the coboundary operator associated to $A^*$,
$$d_* : \Gamma (A) \to \bigwedge\nolimits^2 \Gamma (A),$$
satisfies the cocycle condition (on $\Gamma (A)$)
\begin{equation}\label{Eq:Bialgebroid}
d_* [x, y] = [d_*x, y] + [x, d_*y],
\end{equation}
i.e. if $d_*$ is a derivation of the Schouten algebra $(\Gamma (\bigwedge\nolimits^\bullet A), [\cdot, \cdot])$, we say that the pair $(A, A^*)$ is a \emph{Lie bialgebroid}. Here we use the coboundary operator $ d_*: \Gamma(\bigwedge\nolimits^{k}A) \rightarrow \Gamma(\bigwedge\nolimits^{k+1}A)$ and the Schouten bracket 
given by
\begin{align*}
             d \alpha (x_1, x_2, \ldots, x_{k+1})&= \sum_{i}(-1)^{i-1}\mathrm{a}(x_i)\alpha(x_1,\ldots, \hat{x}_i,\ldots, x_{k+1} )\\
            &+ \sum_{1<i<j<k+1}(-1)^{i+j}\alpha([x_i,x_j], x_1, \ldots, \hat{x}_i,\ldots, \hat{x}_j, \ldots, x_{k+1}).
\end{align*}
and
 \begin{multline*}
            [x_1\wedge \cdots \wedge x_k, y_1 \wedge \cdots \wedge y_l]=\\
            \sum_{i,j} (-1)^{i+j}[x_i,y_j]\wedge x_1 \wedge \cdots \wedge \hat{x}_i \wedge \cdots \wedge x_k \wedge y_1 \wedge \cdots \wedge y_j \wedge \cdots \wedge y_m,        
\end{multline*}
for $x_1 \wedge \cdots \wedge x_k \in \Gamma(\Lambda^kL)$, $y_1 \wedge \cdots \wedge y_m \in \Gamma(\Lambda^mL)$,     and $[x,f] = -[f,x]= \mathrm{a}(x)f$ for $f\in C^{\infty}(M)$,
respectively. \\ \\
It has been shown \cite{Kosmann}\cite{M}\cite{MX} that pairs defining Lie bialgebroids are naturally self-dual: $(A, A^*)$ is a Lie bialgebroid if and only if $(A^*, A)$ is a Lie bialgebroid, and those conditions are equivalent to the coboundary operator associated to $A$ to satisfy the cocycle condition (on $\Gamma (A^*)$), i.e. to $d  : \Gamma (A^*) \to \bigwedge\nolimits^2 \Gamma (A^*)$ to be derivation of the Schouten algebra $(\Gamma (\bigwedge\nolimits^\bullet A^*), [\cdot, \cdot]_*)$ associated to $A^*$. Moreover, the bialgebroid condition can also be written in terms of Lie algebroid morphisms and Poisson maps, namely it is equivalent to the Poisson bundle map $\pi_A^{\#}: T^*A \to TA$ to induce a Lie algebroid morphism with respect to the Lie algebroid structure on $T^*A \to A^*$, see theorem 6.2 in \cite{MX}.\\
\\
Using the coboundary operators $d$ and $d_*$, associated to a dual pair $(A, A^*)$ defining a Lie bialgebroid, it is possible to associate to the space of smooth maps on the base manifold $M$ a Poisson structure \cite{Kosmann}:
\begin{proposition}
Let $(A, A^*)$ be a Lie bialgebroid over a smooth manifold $M$. Then, for any $f, g \in C^\infty (M) = \Gamma (\bigwedge\nolimits^0 A) = \Gamma (\bigwedge\nolimits^0 A^*)$,
\begin{equation}\label{Eq:Poisson-Bialgebroid}
\{f,g\}_{(A, A^*)}= \langle df, d_*g\rangle,
\end{equation}
where $\langle \cdot, \cdot \rangle$ denotes a duality pairing between $A$ and $A^*$, defines a Poisson bracket structure satisfying
$$ d\{f,g\}_{(A, A^*)}= [df, dg]_*.$$
\end{proposition}
\begin{example} In the case of a Poisson manifold $(M, \Pi)$, both $TM$ and $T^*M$ define Lie algebroids and, moreover, they are in duality. The (graded) Jacobi identity makes of the coboundary operator $d_* = [\Pi, \cdot]$  a derivation, so that $(TM, T^*M)$ is a Lie bialgebroid over $M$, sometimes called the standard Lie bialgebroid of  $(M, \Pi)$. The Poisson bracket associated to $f, g \in C^\infty (M)$ is, as expected,
$$ \{f,g\}_{(TM, T^*M)}= \Pi(df,dg).$$
\end{example}
In general, given a Lie bialgebroid $(A, A^*)$ over a smooth manifold $M$ and a function $f \in C^\infty(M)$, a \emph{hamiltonian field} will be any element $X_f \in \Gamma(A)$ satisfying
\begin{equation}\label{Eq:HamVecField}
X_f + d_*f =0.
\end{equation}
In our last example, this corresponds exactly to the usual Hamiltonian vector fields of Poisson geometry.

\section{T-dual pairs and Lie bialgebroids}\label{T-duality}
Given $T$-dual pairs $(E, H)$ and $(\hat{E}, \hat{H})$ (of principal circle bundles with integral $H$-fluxes), we will denote by $A({E})$ and $A(\hat{E})$, respectively, their associated Atiyah algebroids. Recall that, in the topological setting, for $(E, H)$ there exists a closed 2-form $\hat{F} \in \Omega^2(M, \mathfrak{t}^*)$ with integral periods such that,  for all $X \in \mathfrak{t}$,
    \begin{equation}\label{eq:FCurv}
       i_X H = \hat{F}(X),
    \end{equation}
where $\hat{F}(X)$ means the natural pairing (with values in $\Omega^2(M)$) between $\mathfrak{t}$ and $\Omega^2(M,\mathfrak{t}^*)$. The 2-forms defined by (\ref{eq:FCurv}) corresponding to $T$-dual pairs have a natural geometric interpretation as curvatures of the principal $\mathbb{S}^1$-bundles $E$ and $\hat{E}$ and $\Phi = - A \wedge \hat{A}$, where $A$ and $\hat{A}$ denote connections on the $\mathbb{S}^1$-bundles $E$ and $\hat{E}$, satisfying $F=dA$ and $\hat{F}=d \hat{A}$, respectively (see \cite{BHM1} \cite{BHM2} and references therein).\\
\\
The $\mathbb{S}^1$-bundles $E$ and $\hat{E}$ define contact structures on the total spaces, and  $F$ and $\hat{F}$ define symplectic forms on the base space $M$ (see \cite{BW}), we will prove that the Atiyah algebroids $L= A({E})$ and $\hat{L}= A(\hat{E})$ over $M$ define a Lie bi-algebroid, i.e. they are dual 
with anchor $a$ and bracket $[\cdot, \cdot]$, $a_*$ and bracket $[\cdot, \cdot]_*$, respectively, and $d_L$ is a derivation of the Schouten bracket of $\hat{L}$:
$\extd[\alpha, \beta]_* = [\extd \alpha, \beta]_* + [\alpha, \extd\beta]_*$ for all $\alpha, \beta \in \Gamma(\hat{L})$.
The duality pairing giving the identification of  $\hat{L}$ with $L^*$ is the defined by 
\begin{equation}\label{Eq:Duality}
 \llangle X, Y \rrangle  :=\Phi(\Tilde{X}^v,\Tilde{Y}^v) + F(X,Y),
\end{equation}
where $\Tilde{X}$ and $\Tilde{Y}$ denote horizontal lifts \emph{to the correspondence space} $E\times_M\hat{E}$ of the vector fields on each fibration. This pairing is non-degenerate, as follows from the fact that $F$ is non-degenerate and acts only on the horizontal part of vector fields and since, on the other hand, $(E,H)$ and $(\hat{E}, \hat{H})$ are $T$-dual so that the form $\Phi$ is non-degenerate on the vertical part of the corresponding extensions, by definition. 
\\ \\
On the Atiyah algebroid  $\hat{L}= A(\hat{E})$ over $M$ we will consider the local frame $\left\{ \partial_\phi, X_{x^1}, \ldots,  X_{x^n}\right\}$, where $(x_1, \dots, x_n)$ denotes local coordinates on a point in the base and the vector fields $X_{x^j}$, for $j= 1, \dots , n$, are the Hamiltonian vector fields associated to the local coordinates  with respect to the Dirac structure $\mathbb{L}_{F} = \left\{ X + i_X F \mid X \in \mathfrak{X}(M)\right\}$ on $\mathbb{T}M$ associated to the curvature $F$ (the contact structure on $\hat{E}$ can be also used, equivalently, see e.g. \cite{Cardona}\cite{C}. The same can be done for $L= A({E})$ as well). Using the pairing (\ref{Eq:Duality}) we can easily compute the differential $\extd$ (of a smooth function on $M$ or a section of $L^*$) as the next lemma shows.
\begin{lemma} Let $f \in C^\infty (M)$, then 
    \begin{enumerate}
        \item[(a)] $\extd f = X_f$ and $\extd X_f = 0$.
        \item[(b)] $\extd \partial_{\phi} = 0$.
    \end{enumerate}
\end{lemma}
\begin{proof}
    For $(\mathrm{a})$ let $X \in \Gamma(L)$, then 
    \begin{equation*}
        \extd f(X) = a(X) f = df(a(X)) = F(X_f,X) = \llangle X_f, X \rrangle.
    \end{equation*}
    The second part is just the fact that  $\extd^2 = 0$. To prove $(\mathrm{b})$ let $X_1, X_2 \in \mathfrak{X}(M)$, then 
    \begin{align*}
        \extd \partial_\phi (X_1,X_2) = \llangle X_1 \partial_\phi , X_2 \rrangle - \llangle X_2 \partial_\phi , X_1 \rrangle - \llangle \partial_\phi ,([X_1,X_2]) \rrangle= 0.
    \end{align*}
    and
    \begin{align*}
        \extd \partial_\phi(f\partial_\theta, X_1) &= - \llangle X_1\partial_\phi , f\partial_\theta \rrangle - \llangle \partial_\phi , [f\partial_\theta, X_1] \rrangle\\
        &= \llangle -X_1f- \partial_\phi , -X_1f\partial_\theta \rrangle = 0.
    \end{align*}
\end{proof}
Using the preceding lemma it follows that (\ref{Eq:Bialgebroid}) is true for the vector fields of the frame described before. More explicitly, the following equations hold for $j, k= 1, \dots , n$:
\begin{align}
    \extd [X_{x^j}, X_{x^k}]_* &= [\extd X_{x^j}, X_{x^k}]_* + [ X_{x^j}, \extd X_{x^k}]_*,\label{Eq:5.11}\\
    \extd [\partial_\phi, X_{x^j}]_* &= [\extd \partial_\phi, X_{x^j}]_* + [\partial_\phi, \extd X_{x^j}]_*\label{Eq:5.12}. 
\end{align}

\begin{theorem}\label{Theorem}
    The pair $(A(E), A(\hat{E}))$ defines a Lie bi-algebroid.
\end{theorem}

\begin{proof}
In order to show that the pair $(L= A(E), \hat{L} = A(\hat{E}))$ defines a Lie bialgebroid, we need to show that (\ref{Eq:Bialgebroid}) holds for arbitrary sections of $L^*$. By linearity we only need to show that they remain true if we multiply the $X_{x^j}$'s and $\partial_\phi$'s in the equations above by smooth function in $C^\infty(M)$. To this end, we will use that the Poisson structure $\pi$ defined by $F$, satisfy the well-known relations 
\begin{equation*}
    X_{fg} = fX_g + gX_f, \hspace{1cm}  X_f (g) = \{f,g\}
\end{equation*}
and we will write $\{ x^j, x^k\} = \pi^{jk}$. 
Thus, on the left hand side of (\ref{Eq:5.11}) we have
\begin{align*}
    \extd [fX_{x^j}, gX_{x^k}] &=\extd ((fX_{x^j}g)X_{x^k} - (gX_{x^k}f)X_{x^j} + fg[X_{x^j}, X_{x^k}])\\
    &=\extd ((fX_{x^j}g)X_{x^k} - (gX_{x^k}f)X_{x^j} + fgX_{\pi^{jk}})\\
    &= X_{fX_{x^j}g} \wedge X_{x^k} - X_{gX_{x^k}f}\wedge X_{x^j} + X_{fg}\wedge X_{\pi^{jk}},
\end{align*}
while on the right hand side we obtain
\begin{align*}
    [\extd fX_{x^j}, gX_{x^k}]_* &+ [ fX_{x^j}, \extd gX_{x^k}]_*\\
   &= [X_f\wedge X_{x^j}, gX_{x^k}]_* + [fX_{x^j}, X_g\wedge X_{x^k}]_* \\
      &= [X_f, gX_{x^k}]_*\wedge X_{x^j} - [X_{x^j},gX_{x^k}]_*\wedge X_f\\
      & \phantom{=\ } + [fX_{x^j},X_g]_*\wedge X_{x^k}- [fX_{x^j}, X_{x^k}]_*\wedge X_g\\
      &=\{ f,g\} X_{x^k}\wedge X_{x^j} + gX_{\{f,x^k\}}\wedge X_{x^j} - \{ x^j,g \}X_{x^k}\wedge X_f\\
      & \phantom{=\ } - gX_{\pi^{jk}}\wedge X_f - \{ g,f \}X_{x^j}\wedge X_{x^k} + fX_{\{x^j,g\}}\wedge X_{x^k}\\
      & \phantom{=\ } + \{ x^k,f \}X_{x^j}\wedge X_g - fX_{\pi^{jk}}\wedge X_g\\
      &= \{ x^j,g \}X_f\wedge X_{x^k} + fX_{\{x^j,g\}}\wedge X_{x^k}
      -gX_{\{x^k,f\}}\wedge X_{x^j} \\&\phantom{=\ } -
      \{ x^k,f\} X_g\wedge X_{x^j}
        + gX_f\wedge X_{\pi^{jk}} + fX_g\wedge X_{\pi^{jk}}
\end{align*}
showing that the equality remains true. For the case of (\ref{Eq:5.12}) the computation on the left hand side yields
\begin{align*}
    \extd [f\partial_\phi, gX_{x^j}]_* = \extd (-gX_{x^j}f\partial_\phi) = - X_{gX_{x^j}f}\wedge \partial_\phi,
\end{align*}
while on the right-hand side we have
\begin{multline*}
    [\extd (f\partial_\phi), gX_{x^j}]_* + [f\partial_\phi, \extd (gX_{x^j})]_*
    \\
    \begin{aligned}
    &= [X_f\wedge \partial_\phi, gX_{x^j}]_* + [f\partial_\phi, X_g\wedge X_{x^j}]_* \\
    &= [X_f, gX_{x^j}]_*\wedge \partial_\phi - [\partial_\phi,gX_{x^j}]_*\wedge X_f + [f\partial_\phi, X_g]\wedge X_{x^j}- [f\partial_\phi, X_{x^j}]_*\wedge X_g\\
    &= \{f,g\}X_{x^j}\wedge \partial_\phi - gX_{\{x^j,f \}}\wedge \partial_\phi - \{ g,f\} \partial_\phi\wedge X_{x^j} - \{x^j,f\} X_g\wedge X_{x^j}\\
    &= - X_{g\{x^j,f\}}\wedge \partial_\phi.
    \end{aligned}
\end{multline*}
\end{proof}

This result gives us new examples of Lie bi-algebroids and, therefore, of Courant algebroids. Indeed, as shown in \cite{LWX}, any Lie bi-algebroid $(A, A^*)$ gives rise to a canonical Courant algebroid over the same base manifold, namely $\mathbb{E} = A \oplus A^*$. The Courant algebroid  structure of $\mathbb{E}$ comes from the following natural combination of the operations defined on each Lie algebroid: for $x + \alpha, y + \beta \in \Gamma(\mathbb{E})$, the pairing
        \begin{equation*}
            \langle x + \alpha, y + \beta\rangle = \beta(x) + \alpha(y), 
        \end{equation*}
the anchor 
    \begin{equation*}
        \mathbf{a}(x + \alpha) = a(x) + a_*(\alpha) 
    \end{equation*}
and the bracket
    \begin{equation}\label{Eq:4.13}
      \bo x+ \alpha, y + \beta \bc = [x,y] + \Lie_\alpha y - i_{\beta}\extdd x
        + [\alpha, \beta]_* + \Lie_x\beta - i_y\extd\alpha
    \end{equation}
on sections of $\mathbb{E}$, with the natural notations.

\begin{corollary}
    The vector bundle $\mathbb{A}= A(E) \oplus A(\hat{E})$ admits a structure of Courant algebroid over $M$.
\end{corollary}

\begin{example}
    Over $\mathbb{S}^2$ consider the Lie algebroids $T\mathbb{S}^2\oplus\mathfrak{t}$ and $A(\mathbb{S}^3)$, i.e. the Atiyah algebroids of the $T$-dual pairs $(\mathbb{S}^2\times \mathbb{S}^1, H)$ and $(\mathbb{S}^3,0)$ respectively. Then they define a Lie bi-algebroid:
		$$
\begin{tikzcd}
T\mathbb{S}^2\oplus\mathfrak{t}  \arrow[r, "\mathrm{a}"] \arrow[rd] &  T\mathbb{S}^2  \arrow[d]  & A(\mathbb{S}^3) \arrow[l,swap, "\mathrm{a}^*"] \arrow[dl]&  \\
& \mathbb{S}^2.  
\end{tikzcd}
$$
\end{example}
\begin{example} Let $M$ be a compact manifold and let $E$ be a principal $\mathbb{S}^1$-bundle over $M$, with $F$ a representative for $c_1(E)$. Consider the pair $(E, H)$, where $H\in \Omega^3(M)$. Then, for any $X\in \mathfrak{t}$ we obtain $i_XH = 0$ and, hence, if we take $\hat{F} = 0$ we obtain as $T$-dual $\hat{E}$ the trivial bundle $\mathbb{S}^1\times M$. Now endow such bundle with the connection $\hat{A} = d\varphi \otimes \partial_{\varphi}$ and observe that, with the notations used in the introduction, $\Omega = -H$ and in consequence
        \begin{equation*}
            \hat{H} = F\sqcap \hat{A} + H = F\wedge d\varphi + H = d\varphi \wedge F + H.
        \end{equation*}
Hence, the pair $(E, H)$ is $T$-dual to $(\mathbb{S}^1 \times M, d\varphi \wedge F + H)$ and in particular, $(E,0)$ has as $T$-dual $(\mathbb{S}^1\times M, d\varphi \wedge F)$. Thus, considering the associated Atiyah algebroids, $(A(E), A(\mathbb{S}^1\times M))$ defines a Lie bialgebroid.
    \end{example}

\begin{remark}
In \cite{CG} it was shown that $T$-duality, in this context, can be understood as an isomorphism of Courant algebroids   $\mathbb{A}_E^{(H)}$ and $\mathbb{A}_{\hat{E}}^{(\hat{H})}$ which, locally, can be identified with the Courant algebroid
    \begin{equation*}
        T M \oplus T^*M \oplus \mathfrak{t}\oplus \mathfrak{t}^*.
    \end{equation*}
It is shown there that, as expected, $T$-duality is just a permutation of the $\mathfrak{t}$ and $\mathfrak{t}^*$ factors.
\end{remark}

\vspace{0.5cm}

{\bf Acknowledgements.}  The authors are grateful to Henrique
Bursztyn, Michel Cahen, Simone Gutt, Yoshiaki Maeda and Bernardo Uribe 
 for many stimulating discussions on the geometry of Poisson manifolds. This research has been
supported by the \emph{Faculty of Sciences} of the Universidad de los Andes.

\end{document}